\pdfoutput=1
\documentclass[reqno]{amsart}
\usepackage{hyperref, bbm, amssymb, mathtools, url, enumerate}
\usepackage{amsmath}
\usepackage{tikz}
\usepackage{tikz-cd}
\usepackage[scr]{rsfso}
\usepackage{bm}
\usepackage{amssymb}
\usepackage{bbm}
\usetikzlibrary{matrix,arrows,decorations.pathmorphing, decorations.markings, cd}
\usepackage[capitalise]{cleveref}
\usepackage{rotate}
\usepackage[alphabetic]{amsrefs}
\usepackage{microtype}
\usepackage[bb=stixtwo]{mathalpha}
\usepackage[symbol,marginal]{footmisc}

\DeclareMathAlphabet{\mymathbb}{U}{BOONDOX-ds}{m}{n}
\def\twocell[#1]{\arrow[#1, dash, phantom, "\Rightarrow"{scale=1.125, yshift=-.4pt, description, allow upside down, sloped, inner sep=0pt}]}

\tikzset{curve/.style={settings={#1},to path={(\tikztostart)
			.. controls ($(\tikztostart)!\pv{pos}!(\tikztotarget)!\pv{height}!270:(\tikztotarget)$)
			and ($(\tikztostart)!1-\pv{pos}!(\tikztotarget)!\pv{height}!270:(\tikztotarget)$)
			.. (\tikztotarget)\tikztonodes}},
	settings/.code={\tikzset{quiver/.cd,#1}
		\def\pv##1{\pgfkeysvalueof{/tikz/quiver/##1}}},
	quiver/.cd,pos/.initial=0.35,height/.initial=0}

\pgfdeclarelayer{bg}
\pgfsetlayers{bg,main}
\usetikzlibrary{calc}

\newtheorem{theorem}{Theorem}[section]

\newtheorem{corollary}[theorem]{Corollary}
\newtheorem{lemma}[theorem]{Lemma}
\newtheorem{proposition}[theorem]{Proposition}

\newtheorem*{introthm*}{Theorem}

\theoremstyle{definition}
\newtheorem*{claim*}{Claim}
\newtheorem{construction}[theorem]{Construction}

\newtheorem*{remark*}{Remark}

\newcommand{\Cc}{{\mathcal{C}}}

\renewcommand{\phi}{\varphi}
\renewcommand{\epsilon}{\varepsilon}

\DeclareMathOperator{\Spc}{Spc}
\newcommand{\Cat}{\smash{\textup{{C\kern-.75ptat}}}\vphantom{t}}
\newcommand{\CAT}{\smash{\textup{\textsc{C\kern-.75ptat}}}\vphantom{t}}

\DeclareMathOperator{\Fun}{Fun}
\DeclareMathOperator{\PSh}{PSh}

\newcommand{\PSH}{\textsc{PSh}}
\newcommand{\FUN}{\textsc{Fun}}

\newcommand{\op}{{\textup{op}}}

\newcommand{\id}{\textup{id}}

\newcommand{\co}{\textup{co}}

\newcommand{\Fin}{\hbox{$\mathcal F\kern-1.7pt\textit{in}$}}

\newcommand{\incl}{{\mathop{\textup{incl}}}}

\newcommand{\iso}{\xrightarrow{\;\smash{\raisebox{-0.25ex}{\ensuremath{\scriptstyle\sim}}}\;}}

\tikzcdset{pullback/.style = {phantom, "\lrcorner"{very near start}}}

\title[Revisiting $(\infty,2)$-naturality of the Yoneda embedding]{Revisiting $\bm{(\infty,2)}$-naturality of\\ the Yoneda embedding}
\author{Tobias Lenz}
\address{Mathematisches Institut, Rheinische Friedrich-Wilhelms-Universität Bonn, Endenicher Allee 60, 53115 Bonn, Germany}

\begin{document}
\begin{abstract}
	We show that the Yoneda embedding `is' $(\infty,2)$-natural with respect to the functoriality of presheaves via left Kan extension, refining the $(\infty,1)$-categorical result proven independently by Haugseng--Hebestreit--Linskens--Nuiten and Ramzi, and answering a question of Ben-Moshe.

	As the key technical ingredient, we show that the identity functor of the $(\infty,1)$-category of $(\infty,1)$-categories admits only one enhancement to an $(\infty,2)$-functor (namely, the identity functor).
\end{abstract}

\thanks{\\[-1ex]The author is an associate member of the Hausdorff Center for Mathematics at the University of Bonn (DFG GZ 2047/1, project ID 390685813).\\[-1ex]}
\maketitle 

\section*{Introduction}
Even more so than its counterpart in classical category theory, the \emph{Yoneda embedding} is an indispensable tool in higher category theory, for example often allowing us to reduce questions about general $\infty$-categories to the case of (co)complete ones.

In many applications, one needs to be able to perform this reduction `in families,' for which one needs the Yoneda embedding to be \emph{natural} in a suitable sense. It is indeed not hard to show using the universal property of the Yoneda embedding as the free cocompletion \cite{lurie2009HTT} that one can extend the assignment $\Cc\mapsto\PSh(\Cc)$ to a functor $\PSh^\text{free}\colon\Cat_\infty\to\widehat{\Cat}_\infty$  in such a way that the Yoneda embedding becomes natural in an essentially tautological way; unfortunately though, the functor $\PSh^\text{free}$ remains somewhat mysterious this way. On the other hand, the $\infty$-categories of presheaves naturally come to us as a contravariant functor (with functoriality via restriction), and passing to left adjoints produces another covariant functor $\PSh^\text{LKE}\colon\Cat_\infty\to\widehat{\Cat}_\infty$. However, while it is not hard to check that $\PSh^\text{LKE}$ and $\PSh^\text{free}$ agree on the level of homotopy categories, the question of whether they agree fully coherently, or equivalently whether the Yoneda embeddings assemble into a natural transformation $\incl\to\PSh^\text{LKE}$ from the inclusion, turned out to be surprisingly subtle, and an affirmative answer was only given much more recently by Haugseng--Hebestreit--Linskens--Nuiten \cite{HHLN2022TwoVariable} and Ramzi \cite{ramzi-yoneda}; thanks to this, there is no need to distinguish between $\PSh^\text{free}$ and $\PSh^\text{LKE}$ anymore, and so we will denote both simply by $\PSh$ from now on.

This is however still not the strongest result one could hope for: the $\infty$-category $\Cat_\infty$ canonically upgrades to an $(\infty,2)$-category\footnote{Throughout, we will distinguish $(\infty,2)$-categorical refinements of $(\infty,1)$-categories and functors by using small caps.} $\CAT_\infty$ (with mapping objects the usual functor categories), and one would therefore like to have suitable $(\infty,2)$-naturality statements for the Yoneda embedding. The first such result was recently obtained by Ben-Moshe \cite{benmoshe-Yoneda}, who showed that one can use the universal property of presheaves to upgrade the functor $\PSh$ to an $(\infty,2)$-functor $\PSH^\text{free}$ such that the Yoneda embeddings can be uniquely assembled into an $(\infty,2)$-natural transformation. However, the question how this $(\infty,2)$-enhancement relates to the $(\infty,2)$-functor $\PSH^\text{LKE}\colon\CAT_\infty\to\widehat{\CAT}_\infty$ obtained from the natural contravariant $(\infty,2)$-functoriality of presheaf categories (see Construction~\ref{constr:PSHLKE}) was explicitly left open by him. The purpose of this short note is to settle this last remaining question:

\begin{introthm*}[See Theorem~\ref{thm:mainthm}]
	The Yoneda embeddings assemble into an $(\infty,2)$-natural transformation from the inclusion $\CAT_\infty\hookrightarrow\widehat{\CAT}_\infty$ to $\PSH^\textup{LKE}$.
\end{introthm*}

In fact, we will prove a stronger result: no matter how one upgrades the functor $\PSh$ to an $(\infty,2)$-functor, the Yoneda embeddings will assemble into an $(\infty,2)$-natural transformation in a unique way. In particular, all such $(\infty,2)$-functorial enhancements are canonically equivalent.

\subsection*{Strategy and outline} Let us briefly explain the basic idea of the proof, for which we let $\mathbb P\colon\CAT_\infty\to\widehat{\CAT}_\infty$ be any $(\infty,2)$-functorial enhancement of the presheaf functor. We consider the sub-2-functor $\mathbb Q\subset\mathbb P$ given at each $\Cc\in\CAT_\infty$ by the Yoneda image. The Yoneda embeddings then define an equivalence between the underlying $(\infty,1)$-functor of $\mathbb Q$ and the inclusion $\Cat_\infty\hookrightarrow\widehat{\Cat}_\infty$, so it will be enough to show that the latter admits a unique enhancement to an $(\infty,2)$-functor, or equivalently that this holds for the identity of $\Cat_\infty$.

We will prove this rigidity statement in Section~\ref{sec:rigidity}, where we also use this to prove a uniqueness result for the $(\infty,1)$-natural versions of the Yoneda embedding from \cite{HHLN2022TwoVariable} and \cite{ramzi-yoneda}. We then prove our main result in Section~\ref{sec:main} following the argument outlined above.  

\subsection*{Conventions} Below, we will refer to $(\infty,1)$-categories simply as \emph{1-categories} (or \emph{categories} for short), while $(\infty,2)$-categories will be called \emph{2-categories}. As mentioned before, we will distinguish 2-categorical refinements of 1-categorical constructions by using small caps; in particular, $\Cat\subset\widehat{\Cat}$ will denote the 1-category of small 1-categories and the large 1-category of 1-categories, respectively, while the corresponding 2-categories will be denoted $\CAT\subset\widehat{\CAT}$. 

\subsection*{Acknowledgements} I would like to thank Shay Ben-Moshe, Bastiaan Cnossen, Maxime Ramzi, as well as the anonymous referee for feedback on earlier versions of this article.

\section{Rigidity of the identity functor}\label{sec:rigidity}
In this section we will prove our main technical result, showing that the identity of $\Cat$ admits only one lift to a 2-functor $\CAT\to\CAT$. The proof will rely on understanding both the endomorphisms of the 1-functor $\id_{\Cat}$ as well as those of the $2$-functor $\id_{\CAT}$; we begin with the former:

\begin{lemma}\label{lemma:end-id}
	The space of endomorphisms of $\id_{\Cat}\in\Fun(\Cat,\Cat)$ is contractible.
\end{lemma}

That the space of \emph{automorphisms} of the identity is contractible is asserted as part of \cite[Th\'eor\`em 6.3]{toen} or \cite[Theorem~7.3]{bsp-unicity}, so in that sense the mathematical content of the lemma is that every endomorphism of the identity is in fact invertible. The above full statement isn't any harder to prove, however:

\begin{proof}
	By the complete Segal space model, the inclusion $\Delta\hookrightarrow\Cat$ extends to a Bousfield localization $L\colon\PSh(\Delta)\to\Cat$. Precomposing with $L$ and appealing to the universal property of presheaves we then obtain a fully faithful functor
	\[
		\Fun^\textup{L}(\Cat,\Cat)\lhook\joinrel\xrightarrow{\;L^*\;}\Fun^\textup{L}(\PSh(\Delta),\Cat)\iso\Fun(\Delta,\Cat)
	\]
	given by restriction along the fully faithful inclusion $\Delta\hookrightarrow\Cat$. It will therefore suffice to show that the only natural transformation $\sigma\colon\id_\Delta\to\id_\Delta$ is the identity. Since $\Delta$ is a $(1,1)$-category, this may checked by chasing through elements. We then simply observe that for any $0\le k\le n$ the naturality square
	\[
		\begin{tikzcd}
			{[0]}\arrow[d,"k"']\arrow[r,"\sigma_{[0]}"] & {[0]}\arrow[d,"k"]\\
			{[n]}\arrow[r,"\sigma_{[n]}"'] & {[n]}
		\end{tikzcd}
	\]
	witnesses $\sigma_{[n]}(k)=k$.
\end{proof}

\begin{proposition}\label{prop:only-one-pi0}
	Let $\mathbb F\colon\CAT\to\CAT$ be any 2-functor such that its underlying 1-functor $\iota_1\mathbb F\colon\Cat\to\Cat$ is equivalent to the identity. Then $\mathbb F$ is equivalent to the identity 2-functor.
	\begin{proof}
		Note that the identity of $\CAT$ is corepresented by the terminal category $1$. As $\mathbb F(1)\simeq1$ by assumption, the 2-categorical Yoneda lemma shows that there exists a (unique) 2-natural transformation $\tau\colon\id\to\mathbb F$, see \cite[Proposition~11]{benmoshe-Yoneda}. It remains to show that $\tau$ is an equivalence, which can be checked after passing to underlying 1-functors. The previous lemma then shows that the composite $\id_{\Cat}\to\iota_1\mathbb F\simeq\id_{\Cat}$ is the identity, so $\tau$ is an equivalence by 2-out-of-3.
	\end{proof}
\end{proposition}

In other words, the category of 2-functorial enhancements of $\id_{\Cat}$ is connected. To see that it is even contractible, we will need:

\begin{lemma}
	The category of endomorphisms of $\id\colon\CAT\to\CAT$ in the 2-category $\FUN_2(\CAT,\CAT)$ of 2-functors $\CAT\to\CAT$ is a contractible groupoid.
	\begin{proof}
		By another application of the 2-categorical Yoneda lemma  \cite[Proposition~11]{benmoshe-Yoneda}, we have $\textsc{End}(\id)\simeq\id(1)=1$.
	\end{proof}
\end{lemma}

\begin{corollary}\label{cor:rigid}
	The fiber of the forgetful functor
	\[
		\Fun_2(\CAT,\CAT)\to\Fun(\Cat,\Cat)
	\]
	over the identity is a contractible groupoid. Put more concisely, the identity of $\Cat$ admits only one enhancement to a 2-functor $\CAT\to\CAT$.
	\begin{proof}
		By Lemma~\ref{lemma:end-id}, the functor $1\to\Fun(\Cat,\Cat)$ picking out the identity is fully faithful, so the fiber in question is the full subcategory of $\Fun_2(\CAT,\CAT)$ spanned by those 2-functors $\mathbb F$ with $\iota_1\mathbb F\simeq\id_{\Cat}$. By Proposition~\ref{prop:only-one-pi0}, the identity of $\CAT$ is the essentially unique object of this subcategory, while the previous lemma shows that $\id_{\CAT}$ has trivial endomorphisms.
	\end{proof}
\end{corollary}

In \cite{benmoshe-Yoneda}, Ben-Moshe proves a uniqueness statement for his construction of the 2-natural Yoneda embedding: this is a direct consequence of the observation \cite[Proposition~11]{benmoshe-Yoneda} cited above that the inclusion $\CAT\hookrightarrow\widehat{\CAT}$ is corepresented by the terminal category $1$, so that 2-natural transformations into any 2-functor $\mathbb F$ are classified by $\mathbb F(1)$.

Using the above results, we can deduce the analogous statement for the 1-natural Yoneda embedding, strengthening the uniqueness result from \cite[Corollary~2.6]{ramzi-yoneda}.

\begin{theorem}\label{thm:1-unique}
	The 1-natural transformation $y\colon\incl\to\PSh$ from \cite[Theorem~8.1]{HHLN2022TwoVariable} is the unique transformation for which the functor $1\to\PSh(1)\simeq\Spc$ picks out the terminal object.
	\begin{proof}
		Let $\sigma$ be any such transformation. We first claim that each $\sigma(\Cc)\colon\Cc\to\PSh(\Cc)$ factors through the Yoneda image. For this let $X\in\Cc$ be any object (viewed as a functor $X\colon1\to\Cc$) and consider the naturality square
		\[
			\begin{tikzcd}
				1\arrow[d,"X"']\arrow[r,"\sigma"] & \PSh(1)\arrow[d,"X_!"]\\
				\Cc\arrow[r,"\sigma"'] & \PSh(\Cc)\rlap.
			\end{tikzcd}
		\]
		The top path then sends the unique object to $X_!\sigma(1) \simeq X_!y(1)\simeq y(X)$, so we have $\sigma(X)\simeq y(X)\in\text{ess\,im}(y)$ as claimed.

		We now consider the subfunctor $F$ of $\PSh$ given at $\Cc\in\Cat$ by the Yoneda image. By the above, $\sigma$ factors (necessarily uniquely) through $F$, so it will suffice to show that the space of natural transformations $\incl\to F$ is contractible. As $y$ defines an equivalence $\incl\iso F$, this follows immediately from Lemma~\ref{lemma:end-id}. 
	\end{proof}
\end{theorem}

\section{Naturality of the Yoneda embedding}\label{sec:main}
\begin{proposition}\label{prop:2-nat-most-general}
	Let $\mathbb F\colon\CAT\to\widehat{\CAT}$ be any 2-functor, and let $\sigma\colon\incl\to\iota_1\mathbb F$ be a 1-natural transformation from the inclusion $\Cat\hookrightarrow\widehat{\Cat}$ to the underlying 1-functor of $\mathbb F$. Assume moreover that $\sigma(\Cc)\colon\Cc\to \mathbb F(\Cc)$ is fully faithful for every $\Cc\in\Cat$. Then $\sigma$ admits a unique enhancement to a 2-natural transformation from the inclusion $\textsc{incl}\colon\CAT\hookrightarrow\widehat{\CAT}$ to $\mathbb F$.
	\begin{proof}
		Consider the sub-2-functor $\mathbb G$ of $\mathbb F$ given at $\Cc\in\CAT$ by the essential image of $\sigma\colon\Cc\to\mathbb F(\Cc)$; slightly more formally, we can define this by considering the 1-cocartesian unstraightening of $\mathbb F$ in the sense of \cite{Nuiten2021Straightening} and then passing to the evident full subcategory.
		
		By definition, $\sigma$ defines an equivalence $\incl\iso\iota_1\mathbb G$, and we may equivalently show that this equivalence lifts uniquely to a 2-natural transformation. This however follows at once from Corollary~\ref{cor:rigid} together with full faithfulness of the inclusion $\Fun_2(\CAT,\CAT)\hookrightarrow\Fun_2(\CAT,\widehat\CAT)$.
	\end{proof}
\end{proposition}

Specializing this to the Yoneda embeddings, we see that they `are' 2-natural with respect to \emph{any} 2-functorial enhancement of the 1-functor $\PSh\colon\Cat\to\widehat{\Cat}$. For the proof of our main theorem it then only remains to explain what we actually mean by `the 2-functoriality via left Kan extension.'

\begin{construction}\label{constr:PSHLKE}
	Consider the 2-functor $\PSH^\text{res}\colon\CAT^{\co\,\op}\to\widehat{\CAT}$ defined as $\PSH^\text{res}(\Cc)=\Fun(\Cc^\op,\Spc)$, with functoriality via restriction \cite[Example~4.36]{Nuiten2021Straightening}. Note that this factors through the wide and locally full subcategory $\widehat{\CAT}{}^\text{R}$ of right adjoint functors; we may therefore consider the composite
	\[
		\PSH^\text{LKE}\colon\CAT\xrightarrow{\;(\PSH^\text{res})^{\co\,\op}\;}\big(\widehat{\CAT}{}^\text{R}\big)^{\co\,\op}\iso\widehat{\CAT}{}^\text{L}\hookrightarrow\widehat{\CAT},
	\]
	where the unnamed equivalence is the one from \cite[Theorem~3.1.11]{HHLN-mates}, given informally by the identity on objects and by sending a right adjoint functor to its left adjoint. By definition, the underlying 1-functor $\PSh\colon\Cat\to\widehat{\Cat}$ of $\PSH^\text{LKE}$ agrees with the functor from \cite[Theorem~8.1]{HHLN2022TwoVariable}, encoding the functoriality of presheaves via left Kan extension.
\end{construction}

Proposition~\ref{prop:2-nat-most-general} then immediately implies our main theorem: 

\begin{theorem}\label{thm:mainthm}
	The 1-natural Yoneda embedding $y\colon\incl\to\PSh$ (as uniquely characterized in Theorem~\ref{thm:1-unique}) uniquely promotes to a 2-natural transformation $\textsc{incl}\to\PSH^\textup{LKE}$ of 2-functors ${\CAT\to\widehat\CAT}$.\qed
\end{theorem}

\bibliography{reference}
\end{document}